# Liouville-type theorems for twisted and warped products manifolds

Stepanov Sergey

**Abstract.** In the present paper we prove Liouville-type theorems: non-existence theorems for complete twisted and warped products of Riemannian manifolds which generalize and complement similar results for compact manifolds.

## Introduction

In the present paper we will use a *generalized the Bochner technique* which is based on extensions of the classical theorem of divergence and the maximal principle to complete, non-compact Riemannian manifolds. We will apply this technique to prove Liouville-type theorems for complete, non-compact Riemannian doubly twisted and warped products manifolds. We recall that warped products provide are natural generalization of direct products of Riemannian manifolds. The notion of warped products plays very important roles in differential geometry and in general relativity. For example, we recall that many Riemannian manifolds of nonpositive sectional curvature are obtained by using warped products. On the other hand, we recall also that many basic solutions of the Einstein field equations are warped products. For example, the well known Schwarzschild's vacuum model and Robertson-Walker's expanding model of universe in general relativity are warped products. On the other hand, twisted products provide another natural generalization of direct products. The definition of twisted products extends the notion of warped products in a very natural way.

Theorems which we prove in our paper generalize and complement similar well known results for compact doubly twisted and warped products manifolds. In addition, we consider an application of our results to the theory of projective mappings.

# 1. Doubly twisted products manifolds

Let $(M_1, g_1)$ and $(M_2, g_2)$ be two Riemannian manifolds, $\lambda_i : M_1 \times M_2 \to \mathbb{R}$ be a strictly positive differentiable function and $\pi_i : M_1 \times M_2 \to M_i$ be a canonical or natural projection for an arbitrary $i = 1, 2$. Then the *double-twisted product* $_{\lambda_1} M_1 \times_{\lambda_2} M_2$ of $(M_1, g_1)$ and $(M_2, g_2)$ is the differentiable manifold $M = M_1 \times M_2$ equipped with the Riemannian metric $g = \lambda_1^2 g_1 + \lambda_2^2 g_2$ defined by the following equality $g = (\lambda_1^2 \circ \pi_1) \cdot \pi_1^* g_1 + (\lambda_2^2 \circ \pi_2) \cdot \pi_2^* g_2$, where $*$ denotes the pull-back operator on tensors (see [1]). In this case the functions $\lambda_1$ and $\lambda_2$ are called *twisted* functions. In particular, for $\lambda_1 = \lambda_2 = 1$ we have the *direct product* $(M_1 \times M_2, g_1 \oplus g_2)$.

Leaves $M_1 \times \{y\} = \pi_2^{-1}(y)$ and $\{x\} \times M_2 = \pi_1^{-1}(x)$ are totally umbilical submanifolds of $_{\lambda_1} M_1 \times_{\lambda_2} M_2$ (see [1] and [2]). Therefore, the doubly twisted product $_{\lambda_1} M_1 \times_{\lambda_2} M_2$ carries two orthogonal complementary totally umbilical foliations. If we denote by $\pi_{i*} : TM_1 \times TM_2 \to TM_i$ the natural projection then the vector fields $\xi_V = -\pi_{2*}(\nabla \log \lambda_1)$ and $\xi_H = -\pi_{1*}(\nabla \log \lambda_2)$ are the mean curvature vectors of these foliations (see [1] and [2]). The converse is also true in the following case. Namely, let $(M, g)$ be a simply connected $n$-dimensional Riemannian manifold and $\mathcal{V}$ and $\mathcal{H}$ be orthogonal complementary integrable distributions with totally umbilical integral manifolds. If the maximal integral manifolds of $\mathcal{V}$ are complete and $\dim \mathcal{V} = m \geq 3$, then $(M, g)$ is isometric a doubly twisted product of two maximal integral manifolds of $\mathcal{V}$ and $\mathcal{H}$ (see [3]).

We have proved in [2] that the following relation is satisfied on $_{\lambda_1} M_1 \times_{\lambda_2} M_2$

$$\operatorname{div}(\xi_V + \xi_H) = s_{\text{mix}} - \frac{m-1}{m} \|\xi_V\|^2 - \frac{n-m-1}{n-m} \|\xi_H\|^2. \tag{1.1}$$

In this formula $s_{\text{mix}}$ is the *mixed scalar curvature* of $_{\lambda_1} M_1 \times_{\lambda_2} M_2$ which defined as the scalar function

$$s_{\text{mix}} = \sum_{a=1}^{m} \sum_{\alpha=m+1}^{n} \sec(E_a, E_\alpha)$$

where $\sec(E_a, E_\alpha)$ is the *mixed sectional curvature* in direction of the two-plane $\pi = span\{E_a, E_\alpha\}$ for the local orthonormal basis $\{E_1,...,E_m\}$ of the vertical distribution $\mathcal{V}$ and $\{E_{m+1},...,E_n\}$ is orthonormal basis of horizontal distribution $\mathcal{H}$ at an arbitrary point of $_{\lambda_1}M_1 \times_{\lambda_2} M_2$ (see [4, p. 23]; [5]). It is easy to see that this expression is independent of the chosen bases.

Further, taking $\xi_V = -\pi_{2*}(grad \log \lambda_1)$ and $\xi_H = -\pi_{1*}(grad \log \lambda_2)$ in (1.1), we find that

$$div\left(\pi_{2*}(grad \log \lambda_1) + \pi_{1*}(grad \log \lambda_2)\right) =$$
$$= -s_{mix} + \frac{m-1}{m} g_2^*(grad \log \lambda_1, \nabla \log \lambda_1) + \frac{n-m-1}{n-m} g_1^*(grad \log \lambda_2, \nabla \log \lambda_2) \quad (1.2)$$

where

$$g_2^*(grad \log \lambda_1, grad \log \lambda_1) = g(\pi_{2*}grad \log \lambda_1, \pi_{2*}grad \log \lambda_1),$$
$$g_1^*(grad \log \lambda_2, grad \log \lambda_2) = g(\pi_{1*}grad \log \lambda_2, \pi_{1*}grad \log \lambda_2).$$

Let $_{\lambda_1}M_1 \times_{\lambda_2} M_2$ be a complete, noncompact and simply connected Riemannian manifold such that $\|\pi_{2*}(grad \log \lambda_1) + \pi_{1*}(gard \log \lambda_2)\| \in L^1(M, g)$. We recall here that $div X = 0$ for a smooth vector field $X$ on a connected complete, noncompact and oriented Riemannian manifold $(M, g)$ if $\|X\| \in L^1(M, g)$ and $div X \geq 0$ (or $div X \leq 0$) everywhere on $(M, g)$ (see [6] and [7]). In addition, we also recall that every simply connected manifold $M$ is orientable. Then for $s_{mix} \leq 0$ from (1.2) we conclude that $\lambda_1 = C_1$ and $\lambda_2 = C_2$ for some positive constants $C_1$ and $C_2$. In this case $M_1 \times M_2$ has the metric $g = C_1^2 \cdot g_1^* + C_2^2 \cdot g_2^*$. A doubly twisted product $_{\lambda_1}M_1 \times_{\lambda_2} M_2$ with twisted functions $\lambda_1 = C_1 > 0$ *and* $\lambda_2 = C_2 > 0$ can be considered as a direct product $M_1 \times M_2$ of $(M_1, \bar{g}_1)$ and $(M_2, \bar{g}_2)$ for $\bar{g}_1 = C_1^2 g_1$ and $\bar{g}_1 = C_2^2 g_2$. Summarizing, we formulate the statement which generalizes a theorem on two orthogonal complete totally umbilical foliations on compact Riemannian manifold with non positive mixed scalar curvature that has been proved in [2] and [8].

**Theorem 1**. *Let $(M, g)$ be a doubly twisted product ${}_{\lambda_1}M_1 \times_{\lambda_2} M_2$ of some Riemannian manifolds $(M_1, g_1)$ and $(M_2, g_2)$. If $(M, g)$ is a complete and oriented Riemannian manifold $(M, g)$ with nonpositive mixed scalar curvature $s_{\text{mix}}$ and $\|\pi_{2*}(\text{grad} \log \lambda_1) + \pi_{1*}(\text{gard} \log \lambda_2)\| \in L^1(M, g)$, then the twisted functions $\lambda_1$ and $\lambda_2$ are positive constants $C_1$ and $C_2$, respectively, and therefore, $(M, g)$ is the direct product $M_1 \times M_2$ of $(M_1, \bar{g}_1)$ and $(M_2, \bar{g}_2)$ for $\bar{g}_1 = C_1^2 g_1$ and $\bar{g}_1 = C_2^2 g_2$.*

If ${}_{\lambda_1}M_1 \times_{\lambda_2} M_2$ is a well known *Cartan-Hadamard manifold* (see, for example, [9, p. 90]), i.e. a complete, noncompact, simply connected Riemannian manifold of nonpositive sectional curvature, then the above theorem yields the next

**Corollary 1**. *If a Cartan-Hadamard manifold $(M, g)$ is a doubly twisted product ${}_{\lambda_1}M_1 \times_{\lambda_2} M_2$ such that $\|\pi_{2*}(\text{grad} \log \lambda_1) + \pi_{1*}(\text{gard} \log \lambda_2)\| \in L^1(M, g)$, then the twisted functions $\lambda_1$ and $\lambda_2$ are positive constants and therefore, $(M, g)$ is a direct product $M_1 \times M_2$ of $(M_1, \bar{g}_1)$ and $(M_2, \bar{g}_2)$ for $\bar{g}_1 = C_1^2 g_1$ and $\bar{g}_1 = C_2^2 g_2$.*

## 2. Twisted products and projective submersions of Riemannian manifolds

A doubly twisted product manifold ${}_{\lambda_1}M_1 \times_{\lambda_2} M_2$ is called a *twisted product* if $\lambda_1 = 1$ (see [1]). $(M_1, g_1)$ is called the *base* and the submanifolds $M_1 \times \{y\} = \pi_2^{-1}(y)$ are called *leaves* of the twisted product manifold $M_1 \times_{\lambda_2} M_2$. On the other hand, $(M_2, g_2)$ is called *fiber* and the submanifolds $\{x\} \times M_2 = \pi_1^{-1}(x)$ and called *fibres* of the twisted product manifold $M_1 \times_{\lambda_2} M_2$. In addition, all leaves are totally geodesic submanifolds and all fibres are totally umbilical submanifolds in $M_1 \times_{\lambda_2} M_2$. Converse is also true (see [1]). Namely, let $(M, g)$ be a simply connected semi-Riemannian manifold with two orthogonal and complementary integrable distributions $\mathcal{V}$ and $\mathcal{H}$. Suppose that $\mathcal{V}$ is totally geodesic and with complete integrable manifolds. If $\mathcal{H}$ is totally umbilical then $(M, g)$ is isometric to a twisted product $M_1 \times_{\lambda_2} M_2$. In this case, we can formulate the following

**Theorem 2.** *Let a twisted product $M_1 \times_{\lambda_2} M_2$ be a complete and simply connected Riemannian manifold. If the mixed sectional curvature of $M_1 \times_{\lambda_2} M_2$ is non-negative then it is isometric to a direct product $M_1 \times M_2$.*

**Proof.** Our theorem is a corollary of the main theorem of [10] where was considered two orthogonal complementary integrable distributions $\mathcal{V}$ and $\mathcal{H}$ on a complete Riemannian manifold $(M, g)$. If in addition, $\mathcal{V}$ is totally geodesic and $\sum_{a=1}^{m} sec(E_a, E_\alpha) \geq 0$ for each unite vector field $E_\alpha$ which belongs to $\mathcal{H}$ at each point of $M$, then $\mathcal{H}$ is totally geodesic (see [10]). In this case, by the *de Rham decomposition theorem* (see [11, p. 187]) we conclude that if $(M, g)$ is a simply connected Riemannian manifold then it is isometric to the direct product $(M_1 \times M_2, g_1 \oplus g_2)$ of some Riemannian manifolds $(M_1, g_1)$ and $(M_2, g_2)$ for the Riemannian metric $g_1$ and $g_2$ which induced by $g$ on $M_1$ and $M_2$. This completes the proof of the theorem.

We recall here the definition of *pregeodesic* and *geodesic curves*. Namely, a *pregeodesic curve* is a smooth curve $\gamma: t \in J \subset \mathbb{R} \to \gamma(t) \in M$ on a Riemannian manifold $(M, g)$, which becomes a geodesic curve after a change of parameter. Let us change the parameter along $\gamma$ so that $t$ becomes an *affine parameter*. Then $\nabla_X X = 0$ for $X = d\gamma/dt$, and $\gamma$ is called a *geodesic curve*. By analyzing of the last equation, one can conclude that either $\gamma$ is an immersion, i.e., $d\gamma/dt \neq 0$ for all $t \in J$, or $\gamma(t)$ is a point of $M$.

Let $(M, g)$ and $(\overline{M}, \overline{g})$ be Riemannian manifolds of dimension $n$ and $m$ such that $n > m$. A surjective map $f: (M, g) \to (\overline{M}, \overline{g})$ is a *projective submersion* if it has maximal rank $m$ at any point $x$ of $M$ and if $f(\gamma)$ is a pregeodesic in $(\overline{M}, \overline{g})$ for an arbitrary pregeodesic $\gamma$ in $(M, g)$ (see [12]). In this case, each pregeodesic line $\gamma \subset M$ which is an integral curve of the distribution $Ker f_*$ is mapped into a point $f(\gamma)$ in $\overline{M}$. Note that this fact does not contradict the definition of projective submersion.

We call the submanifolds $f^{-1}(y) \subset M$ for an arbitrary $y \in \overline{M}$ as *fibers*. All fibers of an arbitrary projective submersion are totally geodesic submanifolds (see [12] and [13]). Putting $\mathcal{V}_x = Ker(f_*)_x$, for any $x \in M$, we obtain an integrable vertical distribution $\mathcal{V}$ which corresponds to the totally geodesic foliation of $M$ determined by the fibres of $f$, since each $\mathcal{V}_x = T_x f^{-1}(y)$ coincides with tangent space of $f^{-1}(y)$ at $x$ for $f(x) = y$. Let $\mathcal{H}$ be the complementary distribution of $\mathcal{V}$ determined by the Riemannian metric $g$, i.e. $\mathcal{H}_x = \mathcal{V}_x^\perp$ at each $x \in M$ where $\mathcal{H}_x$ is called the *horizontal space* at $x$
We have proved in [12] and [13] that the horizontal distribution $\mathcal{H}$ is integrable with totally umbilical integral manifolds. So, $(M, g)$ admits two complementary totally geodesic and totally umbilical foliations, whose leaves intersect perpendicularly. If $(M, g)$ is a simply connected Riemannian manifold and the integral manifolds of $\mathcal{H}$ are geodesically complete then $(M, g)$ is isometric to a twisted product $M_1 \times_{\lambda_2} M_2$ such that the maximal integral manifolds of $\mathcal{V}$ and $\mathcal{H}$ correspond to the canonical foliations of the product $M_1 \times M_2$ (see [1]).

The converse is also true in local case. Namely, for an arbitrary twisted product $M_1 \times_{\lambda_2} M_2$ the natural projection $\pi_2 : M_1 \times_{\lambda_2} M_2 \to M_2$ is a locally projective submersion.

We make here two observations before formulating a conclusion. Firstly, if we suppose the geodesic completeness of $(M, g)$, then the integral manifolds of are geodesically complete automatically. Secondary, it well known by the Hopf–Rinow theorem that any connected geodesically complete Riemannian manifold is a complete Riemannian manifold. Then the following statement is a corollary of Theorem 2.

**Corollary 2**. *Let $(M, g)$ be an n-dimensional complete and simply connected Riemannian manifold with non-negative sectional curvature. If $(M, g)$ admits a projective submersion $f : (M, g) \to (\overline{M}, \overline{g})$ onto another m-dimensional ($m < n$) Riemannian manifold $(\overline{M}, \overline{g})$ then it is isometric to a direct product $(M_1 \times M_2, g_1 \oplus g_2)$ of some*

Riemannian manifolds $(M_1, g_2)$ and $(M_1, g_2)$ such that maximal integral manifolds of $\operatorname{Ker} f_*$ and $(\operatorname{Ker} f_*)^\perp$ correspond to the canonical foliations of the product $M_1 \times M_2$.

Another corollary follows from Theorem 1.

**Corollary 3**. *Let $(M, g)$ be an n-dimensional complete and simply connected Riemannian manifold and $f : (M, g) \to (\overline{M}, \overline{g})$ be a projective submersion onto another m-dimensional ($m < n$) Riemannian manifold $(\overline{M}, \overline{g})$, then $(M, g)$ is a twisted product $M_1 \times_{\lambda_2} M_2$ of some Riemannian manifolds $(M_1, g_1)$ and $(M_2, g_2)$. If $M_1 \times_{\lambda_2} M_2$ has a nonpositive mixed scalar curvature $s_{\text{mix}}$ and $\|\pi_{1*}(\operatorname{gard} \log \lambda_2)\| \in L^1(M, g)$, then it isometric to a direct product $(M_1 \times M_2, g_1 \oplus g_2)$ such that maximal integral manifolds of $\operatorname{Ker} f_*$ and $(\operatorname{Ker} f_*)^\perp$ correspond to the canonical foliations of $M_1 \times M_2$.*

## 3. Double warped products and warped products manifolds

A *doubly warped product* manifold $(M, g)$ is a twisted product manifold ${}_{\lambda_1}M_1 \times_{\lambda_2} M_2$ where $\lambda_1 : M_2 \to \mathbb{R}$ and $\lambda_2 : M_1 \to \mathbb{R}$ are positive smooth functions (see [14]). These functions are called *warping* functions. A doubly warped product ${}_{\lambda_1}M_1 \times_{\lambda_2} M_2$ carries two orthogonal complementary totally umbilical foliations with closed mean curvature vectors $\xi_V$ and $\xi_H$ because in this case we have the equalities $\xi_V = -\pi_{2*}(\operatorname{grad} \log \lambda_1) = -\operatorname{grad} \log \lambda_1$ and $\xi_H = -\pi_{1*}(\operatorname{grad} \log \lambda_2) = -\operatorname{grad} \log \lambda_2$, respectively (see also [1]). Then the formula (1.2) can be rewrite in the following form

$$\Delta \log(\lambda_1 \lambda_2) = -s_{\text{mix}} + \frac{m-1}{m} \|\operatorname{grad} \log \lambda_1\|^2 + \frac{n-m-1}{n-m} \|\operatorname{grad} \log \lambda_2\|^2 \qquad (3.1)$$

where $\|\operatorname{grad} \log \lambda_i\|^2 = g(\operatorname{grad} \log \lambda_i, \operatorname{grad} \log \lambda_i)$ for an arbitrary $i = 1, 2$. . Therefore, if $s_{\text{mix}} \le 0$ then the function $\log(\lambda_1 \lambda_2)$ is subharmonic because $\Delta \log(\lambda_1 \lambda_2) \ge 0$. At the same time, we know that on a complete Riemannian manifold $(M, g)$ each subharmonic function $f : M \to \mathbb{R}$ whose gradient has integrable norm on $(M, g)$ must actually be harmonic (see [15, p. 660]). In our case it means that $\Delta \log(\lambda_1 \lambda_2) = 0$. Then

from (3.1) we obtain $\lambda_1 = C_1$ and $\lambda_2 = C_2$ for some positive constants $C_1$ and $C_2$. Therefore, the following theorem holds.

**Theorem 3**. *Let (M, g) be a complete double-warped product $_{\lambda_1}M_1 \times_{\lambda_2} M_2$ of Riemannian manifolds $(M_1, g_1)$ and $(M_2, g_2)$ such that the mixed scalar curvature $s_{\text{mix}}$ is nonpositive. If the gradient of $\log(\lambda_1 \lambda_2)$ has integrable norm (M, g), then $\lambda_1 = C_1$ and $\lambda_2 = C_2$ for some positive constants $C_1$ and $C_2$ and therefore, (M, g) is the direct product of $(M_1, \bar{g}_1)$ and $(M_2, \bar{g}_2)$ for $\bar{g}_1 = C_1 g_1$ and $\bar{g}_1 = C_2 g_2$.*

**Remark**. We argue that our Theorem 3 complements the results of the paper [15] where was proved that if the mixed sectional curvature of a complete doubly warped product manifold $_{\lambda_1}M_1 \times_{\lambda_2} M_2$ is non-negative then the warping functions $\lambda_1$ and $\lambda_2$ are constants.

Twisted products are generalizations of *warped products*. For a warped product $M_1 \times_{\lambda_2} M_2$ the function $\lambda_2$ is a smooth positive function $\lambda_2 : M_1 \to \mathbb{R}$ (see [16] and [17]). In this case, leaves are totally geodesic submanifolds in $M_1 \times_{\lambda_2} M_2$ and fibres are extrinsic spheres. In addition, we recall that a submanifold of a Riemannian manifold is called an *extrinsic sphere* if it is a totally umbilical submanifold with parallel mean curvature vector (see [18]). In [1] was proved the following statement: Let (M, g) be a simply connected semi-Riemannian manifold with two orthogonal and complementary integrable distributions $\mathcal{V}$ and $\mathcal{H}$ such that $\mathcal{V}$ is totally geodesic and with complete leaves and integrable manifolds of $\mathcal{H}$ are extrinsic spheres. Then (M, g) is isometric to a warped product.

Let (M, g) be a warped product $M_1 \times_{\lambda_2} M_2$ of two Riemannian manifolds $(M_1, g_1)$ and $(M_2, g_2)$. The well known curvature identities (see [19, p. 211])

$$\pi_1^* \text{Ric} = \text{Ric}_1 - \frac{n-m}{\lambda_2} \text{Hess}(\lambda_2),$$

$$\pi_2^* \text{Ric} = \text{Ric}_2 - \left(\frac{\Delta_1 \lambda_2}{\lambda_2} - (n-m-1)\frac{\pi_2^* g(\text{grad } \lambda_2, \text{grad } \lambda_2)}{\lambda_2^2}\right) \pi_2^* g,$$

relating the Ricci curvature $Ric$ of $(M, g)$ and the Ricci curvatures $Ric_1$ and $Ric_2$ of $(M_1, g_1)$ and $(M_2, g_2)$, respectively. In these identities $\Delta_1 \lambda_2$ is the Laplacian of $\lambda_2$ for $g_1$ and $\pi_2^* g = \lambda_2^2 g_2$.

From the above identities we obtain two equations

$$\Delta_1 \lambda_2 = \frac{1}{n-m} \lambda_2 \left(s_1 - trace_{g_1}\left(\pi_1^* Ric\right)\right) \tag{3.2}$$

$$\lambda_2 \Delta_1 \lambda_2 = \frac{1}{n-m}\left(s_2 - trace_{g_2}\left(\pi_2^* Ric\right)\right) - (n-m-1)\lambda_2^2 \| grad\, \lambda_2 \|^2 \tag{3.3}$$

where $\| grad\, \lambda_2 \|^2 = g_1(grad\, \lambda_2, grad\, \lambda_2)$, $trace_{g_1}\left(\pi_1^* Ric\right) = \sum_{a=1}^{m} Ric(E_a, E_a)$ and $trace_{g_2}\left(\pi_2^* Ric\right) = \sum_{\alpha=m+1}^{n} Ric(E_\alpha, E_\alpha)$. In first case, if we assume that $s_1 \geq trace_{g_1} \pi_1^* Ric$ then from (3.2) we obtain $\Delta_1 \lambda_2 \geq 0$ and therefore the warped function $\lambda_2 : M_1 \to \mathbb{R}$ is subharmonic. It is well known that Yau showed in [20] that every non-negative $L^p$-subharmonic function on a complete Riemannian manifold must be constant for any $p > 1$. Therefore, if $(M_1, g_1)$ is a complete manifold such that its scalar curvature $s_1 \geq trace_{g_1} Ric$ and $\int_{M_1} \lambda_2^p dV_{g_1} < \infty$ for some $p > 1$ then from (3.2) we conclude that the warped function $\lambda_2$ is constant. At the same time, a warped product $M_1 \times_{\lambda_2} M_2$ with a constant warping function $\lambda_2 = C_2 > 0$ can be considered as a direct product $M_1 \times M_2$ of $(M_1, g_1)$ and $(M_2, \bar{g}_2)$ for $\bar{g}_1 = C_2 g_2$. Summarizing the above arguments we can formulate the following

**Theorem 4**. *Let $(M, g)$ be a warped product $M_1 \times_{\lambda_2} M_2$ of two Riemannian manifolds $(M_1, g_1)$ and $(M_2, g_2)$ such that the base $(M_1, g_1)$ of $M_1 \times_{\lambda_2} M_2$ is a complete manifold and $s_1 \geq trace_{g_1} \pi_1^* Ric$ for the scalar curvature $s_1$ of $(M_1, g_1)$ and for the Ricci tensor $Ric$ of $M_1 \times_{\lambda_2} M_2$. If $\int_{M_1} \lambda_2^p dV_{g_1} < \infty$ for some $p > 1$ then $\lambda_2 = C_2$ for some positive constant $C_2$ and therefore, $(M, g)$ is the direct product $M_1 \times M_2$ of $(M_1, g_1)$ and $(M_2, \bar{g}_2)$ for $\bar{g}_2 = C_2 g_2$.*

In the second case, if we assume that $s_2 \geq trace_{g_2} \pi_2^* Ric$ then from (3.3) we obtain $\Delta_1 \lambda_2 \leq 0$ and therefore the warped function $\lambda_2 : M_1 \to \mathbb{R}$ is superharmonic. In this case, if $(M, g)$ is a complete Riemannian manifold and $\|grad\, \lambda_2\| \in L^1(M, g)$, then superharmonic function $\lambda_2$ is harmonic (see the Lemma in Appendix). Thus, we can formulate an analogue of the previous theorem.

**Theorem 5**. *Let $(M, g)$ be a warped product $M_1 \times_{\lambda_2} M_2$ of two Riemannian manifolds $(M_1, g_1)$ and $(M_2, g_2)$ such that the base $(M_1, g_1)$ of $M_1 \times_{\lambda_2} M_2$ is a complete manifold and $s_2 \leq trace_{g_2} \pi_2^* Ric$ for the scalar curvature $s_2$ of $(M_2, g_2)$ and for the Ricci tensor Ric of $M_1 \times_{\lambda_2} M_2$. If $\|grad\, \lambda_2\| \in L^1(M, g)$ then $\lambda_2 = C_2$ for some positive constant $C_2$ and therefore, $(M, g)$ is the direct product $M_1 \times M_2$ of $(M_1, g_1)$ and $(M_2, \bar{g}_2)$ for $\bar{g}_1 = C_2 g_2$.*

If we assume that $u = \lambda_2^{(n-m+1)/2}$ then the following equation holds (see [21])

$$\frac{4(n-m)}{n-m+1} \Delta_1 u = (s_1 - s)u + s_2\, u^{(n-m-3)/(n-m+1)} \tag{3.4}$$

for $n \geq 3$. In turn, if $s \leq s_1$ and $s_2 \geq 0$ then from the equations (3.4) followed $\Delta_1 u \geq 0$, i.e. $u$ is a subharmonic function on $(M_1, g_1)$. Therefore, if we suppose that $(M_1, g_1)$ is a complete manifold and $u = \lambda_2^{(n-m+1)/2}$ is $L^p$-function on $(M_1, g_1)$ for some $p > 1$ then the warped product function $\lambda_2$ is a positive constant. In particular, for $n = m + 3$ the equation (3.4) can be rewritten in the form $3\Delta_1 u = (s_1 - s)u$. In this case, if $s \leq s_1$ then $u$ is a subharmonic function on $(M_1, g_1)$. Thus we have the following result.

**Theorem 6**. *Let $M_1 \times_{\lambda_2} M_2$ be a warped product such that its base $(M_1, g_1)$ is a complete manifold. If the warping function $\lambda_2^{(n-m+1)/2} \in L^p(M_1, g_1)$ for some $p > 1$ and the scalar curvatures $s$, $s_1$ and $s_2$ of $M_1 \times_{\lambda_2} M_2$ and of its base and fibre, respectively, satisfy one of the two following conditions*:

1. $s \geq s_1$ and $s_2 \leq 0$ for any $n \geq 3$;

2. $s \leq s_1$ for $n - m = 3$, then $\lambda_2$ is a positive constant $C_2$ and therefore, $(M, g)$ is isometric to the direct product $M_1 \times M_2$ of $(M_1, g_1)$ and $(M_2, \bar{g}_2)$ for $\bar{g}_1 = C_2 g_2$.

## 4. Einstein warped product manifolds

Let the warped product $M_1 \times_{\lambda_2} M_2$ be an $n$-dimensional ($n \geq 3$) Einstein manifold, i.e. $Ric = \frac{s}{n} g$ for the constant scalar curvature $s$ of $M_1 \times_{\lambda_2} M_2$. In this case, $(M_2, g_2)$ is an Einstein manifold too (see [16]). It means that for $n - m \geq 3$ we have $Ric_2 = \frac{s_2}{n - m} g_2$ where $Ric_2$ and $s_2$ are the Ricci tensor and the constant scalar curvature of $(M_2, g_2)$, respectively. In addition, the following equation holds

$$\lambda_2 \Delta \lambda_2 = (n - m - 1)\| \operatorname{grad} \lambda_2 \|^2 + \left( \frac{s}{n} \lambda_2^2 - \frac{s_2}{n - m} \right). \tag{4.1}$$

On the other hand, for $Ric = \frac{s}{n} g$ the equation (3.3) can be rewritten in the form

$$\lambda_2 \Delta_1 \lambda_2 = \left( \frac{s_2}{n - m} - \frac{s}{n} \right) - (n - m - 1)\lambda_2^2 \| \operatorname{grad} \lambda_2 \|^2. \tag{4.2}$$

From (4.1) and (4.2) we obtain

$$2\lambda_2 \Delta_1 \lambda_2 = \left( \lambda_2^2 - 1 \right)\left( \frac{s}{n} - (n - m - 1)\| \operatorname{grad} \lambda_2 \|^2 \right). \tag{4.3}$$

and

$$\frac{2 s_2}{n - m} = \left( \lambda_2^2 + 1 \right)\left( \frac{s}{n} + (n - m - 1)\lambda_2^2 \| \operatorname{grad} \lambda_2 \|^2 \right). \tag{4.4}$$

Now let us analyze the equations (4.1) – (4.4). Firstly, if $s > 0$ then from (4.4) we obtain $s_2 > 0$. In this case, from (4.4) we conclude that $\lambda_2^2 \leq 2\frac{n s_2}{(n - m)s} - 1$. Secondary, we consider the equation (4.1). For this case, we recall the well known statement from [20]. Namely, if $f$ is a smooth function defined on a complete Riemannian manifold $(M, g)$ satisfies the equality $f \Delta f \geq 0$, then either $\int_M |f|^p dV_g = 0$ for all $p \neq 1$ or $f = $ constant. It means that if $M_1 \times_{\lambda_2} M_2$ is complete manifold such that $\lambda_2^2 \geq \frac{n s_2}{(n - m)s}$ for

$s > 0$ and $\int_{M_1} \lambda_2^p dV_{g_1} < \infty$ for some $p \neq 1$ then from (4.1) we obtain that $\lambda_2 = \sqrt{\dfrac{ns_2}{(n-m)s}}$ is constant and therefore $(M, g)$ is the direct product $M_1 \times M_2$ of $(M_1, g_1)$ and $(M_2, \bar{g}_2)$ for $\bar{g}_1 = \lambda_2 g_2$. In addition, we note that from the equalities $\lambda_2^2 \geq \dfrac{ns_2}{(n-m)s}$ and $\lambda_2^2 \leq 2\dfrac{ns_2}{(n-m)s} - 1$ follows $\lambda_2^2 \geq \dfrac{ns_2}{(n-m)s} \geq 1$.

Thirdly, we consider the equation (4.3). In this case, if we assume that $s < 0$ and $\int_{M_1} \lambda_2^p dV_{g_1} < \infty$ for $p \neq 1$ and $\lambda_2 \leq 1$, then we conclude that $\Delta \lambda_2 = 0$ and $\lambda_2 = 1$. On the other hand, if we assume that $s < 0$ and $\int_{M_1} \| grad\, \lambda_2 \| dV_{g_1} < \infty$ for $\lambda_2 \geq 1$, then we conclude that $\Delta \lambda_2 = 0$ and $\lambda_2 = 1$ (see the Lemma in Appendix). Summarizing the above arguments we can formulate the following statement.

**Corollary 4**. *Let $M_1 \times_{\lambda_2} M_2$ be an n-dimensional ($n \geq 3$) Einstein warped product of two Riemannian manifolds $(M_1, g_1)$ and $(M_2, g_2)$ such that the m-dimensional ($n - m \geq 3$) base $(M_1, g_1)$ of $M_1 \times_{\lambda_2} M_2$ is complete manifold.*

1) *If $\lambda_2^2 \geq \dfrac{ns_2}{(n-m)s} \geq 1$ for the positive constant scalar curvatures $s_2$ and $s$ of $(M_2, g_2)$ and $M_1 \times_{\lambda_2} M_2$, respectively and $\lambda_2 \in L^p(M_1, g_1)$ for some $p \neq 1$, then $\lambda_2 = \sqrt{\dfrac{ns_2}{(n-m)s}}$ and therefore $(M, g)$ is the direct product $M_1 \times M_2$ of $(M_1, g_1)$ and $(M_2, \bar{g}_2)$ for $\bar{g}_1 = \lambda_2 g_2$.*

2) *If the scalar curvature $s$ of $M_1 \times_{\lambda_2} M_2$ is negative and $\lambda_2 \in L^p(M_1, g_1)$ for some $p \neq 1$ and $\lambda_2 \leq 1$, then $\lambda_2 = 1$ and therefore $(M, g)$ is the direct product $M_1 \times M_2$ of $(M_1, g_1)$ and $(M_2, g_2)$.*

3) *If the scalar curvature s of $M_1 \times_{\lambda_2} M_2$ is negative and $\|grad\, \lambda_2\| \in L^1(M_1, g_1)$ for $\lambda_2 \geq 1$, then $\lambda_2 = 1$ and therefore $(M, g)$ is the direct product $M_1 \times M_2$ of $(M_1, g_1)$ and $(M_2, g_2)$.*

If the warped product $M_1 \times_{\lambda_2} M_2$ is an *n*-dimensional ($n \geq 3$) Einstein manifold, then (3.2) can be rewritten in the form

$$\Delta_1 \lambda_2 = \frac{1}{n-m} \lambda_2 \left( s_1 - \frac{m}{n} s \right), \qquad (4.5)$$

If $(M_1, g_1)$ is a complete manifold such that $s_1 \geq \frac{m}{n} s$ and $\int_{M_1} \lambda_2^p dV_{g_1} < \infty$ for some $p \neq 1$, then from (4.5) we conclude that the warped function $\lambda_2$ is constant and $s_1 = \frac{m}{n} s =$ *constant*. On the other hand, $s_1 \geq \frac{m}{n} s$ and $\int_{M_1} \|grad\, \lambda_2\| dV_{g_1} < \infty$, then from (4.5) we conclude that the warped function $\lambda_2$ is constant and $s_1 = \frac{m}{n} s =$ *constant* (see the Lemma in Appendix). We proved the following

**Corollary 5**. *Let $M_1 \times_{\lambda_2} M_2$ be an n-dimensional ($n \geq 3$) Einstein warped product of two Riemannian manifolds $(M_1, g_1)$ and $(M_2, g_2)$ such that the base $(M_1, g_1)$ of $M_1 \times_{\lambda_2} M_2$ is an m-dimensional complete manifold and $s_1 \geq \frac{m}{n} s$ (resp. $s_1 \leq \frac{m}{n} s$) for the scalar curvature $s_1$ of $(M_1, g_1)$ and for the constant scalar curvature s of $M_1 \times_{\lambda_2} M_2$. If $\lambda_2 \in L^p(M_1, g_1)$ for some $p \neq 1$ (resp. $\|grad\, \lambda_2\| \in L^1(M_1, g_1)$) then $s_1 = \frac{m}{n} s =$ constant and $\lambda_2 = C_2$ for some positive constant $C_2$ and therefore $(M, g)$ is the direct product $M_1 \times M_2$ of $(M_1, g_1)$ and $(M_2, \bar{g}_2)$ for $\bar{g}_1 = C_2 g_2$.*

**Remark**. Corollaries 4 and 5 generalize and complement the main theorem on an Einstein warped product with compact base in [16].

## 5. Appendix.

In conclusion, we consider *superharmonic function* on complete, non-compact Riemannian manifolds and prove the following lemma which is an analogy of the Yau proposition from [15, p. 660]) in which he has argued that on a complete non-compact Riemannian manifold each subharmonic function whose gradient has integrable norm on $(M, g)$ must be harmonic.

**Lemma**. *If $(M, g)$ is a complete Riemannian manifold without boundary, then any superharmonic function $f \in C^2 M$ with $\|\operatorname{grad} f\| \in L^1(M, g)$ is harmonic.*

**Proof**. On the one hand, if we assume that $\varphi = -f$ for any superharmonic function $f \in C^2 M$ then the conditions $\Delta f \leq 0$ and $\|\operatorname{grad} f\| \in L^1(M, g)$ which must be satisfy for the superharmonic function $f$ can be written in the form $\Delta \varphi \geq 0$ and $\|\operatorname{grad} \varphi\| \in L^1(M, g)$. In this case, using the Yau statement we conclude that $\Delta \varphi = 0$ and hence $f = -\varphi$ is a harmonic function.

## Acknowledgements

The work of Sergey Stepanov is supported by RBRF grant 16-01-00053-a (Russia).

## References

[1] Ponge R., Reckziegel H., Twisted products in pseudo-Riemannian geometry, Geom. Dedic., 48:1 (1993), 15-25.

[2] Stepanov S.E., A class of Riemannian almost-product structures, Soviet Mathematics (Izv. VUZ), 33:7 (1989), 51-59.

[3] Kim B.H., Warped products with critical Riemannian metric, Proc. Japan Acad., Ser. A, 71:6 (1995), 117-118.

[4] Falcitelli M., Ianus S., Pastore A.M., Riemannian submersions and related topics, Word Scientific Publishing, Singapore, 2004.

[5] Rocamora A.H., Some geometric consequences of the Weitzenböck formula on Riemannian almost-product manifolds; weak-harmonic distributions, Illinois Journal of Mathematics, 32:4 (1988), 654-671.

[6] Caminha A., Souza P., Camargo F., Complete foliations of space forms by hypersufaces, Bull. Braz. Math. Soc., New Series, 41:3 (2010), 339-353.

[7] Caminha A., The geometry of closed conformal vector fields on Riemannian spaces, Bull. Braz. Math. Soc., New Series, 42:2 (2011), 277-300.

[8] Naveira A.M., Rocamora A.H., A geometrical obstruction to the existence of two totally umbilical complementary foliations in compact manifolds, Differential Geometrical Methods in Mathematical Physics, Lecture Notes in Mathematics 1139 (1985), 263-279.

[9] Pigola S., Rigoli M., Setti A.G., Vanishing and Finiteness Results in Geometric Analysis. A Generalization of the Bochner Technique, Birkhäuser Verlag AG, Berlin (2008).

[10] Brito F., Walczak P.G., Totally geodesic foliations with integrable normal bundles, Bol. Soc. Bras. Mat., 17:1 (1986), 41-46.

[11] Koboyashi S., Nomizu K., Foundations of differential geometry, Volume I, Interscience Publishers, New York, 1963.

[12] Stepanov S.E., On the global theory of projective mappings, Mathematical Notes, 58:1 (1995), 752-756.

[13] Stepanov S.E., Geometry of projective submersions of Riemannian manifolds, Russian Mathematics (Iz. VUZ), 43:9 (1999), 44-50.

[14] Ünal B., Doubly warped products, Differential Geometry and its Applications, 15 (2001), 253-263.

[15] Gutierrez M., Olea B., Semi-Riemannian manifolds with a doubly warped structure, Revista Matematica Iberoamericana, 28:1 (2012), 1-24.

[16] Kim D.-S., Kim Y.H., Compact Einstein warped product spaces with nonpositive scalar curvature, Proceedings of the American Mathematical Society, 131:8 (2003), 2573-2576.

[17] Defever F., Deszcz R., Glogowska M., Goldberg V.V., Verstraelen L., A class of four-dimensional warped products, Demonstatio Mathematica, 35:4 (2002), 853-864.

[18] Kozaki M., Ohkubo T., A characterization of extrinsic spheres in a Riemannian manifold, Tsukuba J. Math., 26:2 (2002), 291-297.

[19] O'Neill B., Semi-Riemannian geometry with applications to relativity, Academic Press, San Diego, 1983.

[20] Yau S.T., Some function-theoretic properties of complete Riemannian manifolds and their applications to geometry, Indiana Univ. Math. J., 25 (1976), 659-670.

[21] Dobarro F., Dozo E.L., Scalar curvature and warped products of Riemannian manifolds, Transactions of the American Mathematical Society, 303:1 (1987), 161-168.